\documentclass[12pt]{article}
\usepackage[affil-it]{authblk}
\usepackage{amsmath} 
\usepackage{amssymb} 
\usepackage{graphicx} 
\usepackage{enumerate}
\usepackage[utf8]{inputenc}
\usepackage{amsthm}
\usepackage{color}

\newtheorem{thm}{Theorem}

\newtheorem{prop}[thm]{Proposition}

\newtheorem{ex}[thm]{Example}

\newcommand{\F}{\mathbb{F}}
\newcommand{\Fq}{\mathbb{F}_q}
\newcommand{\Fp}{\mathbb{F}_p}
\newcommand{\On}{\textbf{O}_{n}}

\begin{document}

\title{Complete mappings and Carlitz rank}

\date{}
\author{Leyla I\c s\i k$^1$,  Alev Topuzo\u glu$^1$, Arne Winterhof$^2$}

\maketitle

\noindent
$^1$ Sabanc\i~ University, Orhanli, 34956 Tuzla, \.{I}stanbul, Turkey\\
E-mail: \{isikleyla,alev\}@sabanciuniv.edu\\

\noindent
$^2$ Johann Radon Institute for Computational and Applied Mathematics\\
Austrian Academy of Sciences, Altenbergerstr.\ 69, 4040 Linz, Austria\\
E-mail: arne.winterhof@oeaw.ac.at

\begin{abstract}
The well-known Chowla and Zassenhaus conjecture, proven by Cohen in 1990, states that 
for any $d\ge 2$ and any prime $p>(d^2-3d+4)^2$ there is no complete mapping polynomial in 
$\Fp[x]$ of degree $d$. 

For arbitrary finite fields $\Fq$,
we give a similar result in terms of the Carlitz rank
of a permutation polynomial  
rather than its degree. 
We prove that if $n<\lfloor q/2\rfloor$, then  
there is no complete
mapping in $\Fq[x]$ of Carlitz rank $n$ of small linearity. We also
determine how far permutation polynomials $f$ of Carlitz rank $n<\lfloor q/2\rfloor$ are from being complete, by studying
value sets of $f+x.$ We provide examples of complete mappings 
if $n=\lfloor q/2\rfloor$, which shows that the above bound 
cannot be improved in general. 
\end{abstract}

\bigskip


\bigskip


\section{Introduction}
For any prime power $q$ let
$\F_q$ be the finite field of $q$ elements. 
A polynomial $f(x) \in\F_q[x]$ is called a {\it permutation
polynomial} if it induces a bijection from~$\F_q$ to~$\F_q$.

A polynomial $f(x)\in\mathbb{F}_q[x]$ is a \textit{complete mapping
polynomial} (or a complete mapping) if both $f(x)$ and $f(x)+x$ are permutation 
polynomials of~$\mathbb{F}_q$. These polynomials were introduced by Mann in 1942,
\cite{Mann1942}. A detailed study of complete mapping polynomials
over finite fields was carried out by Niederreiter and Robinson
(1982, \cite{NR}). 
Complete mappings are pertinent to the construction of mutually orthogonal Latin squares, which can be used for the design of agricultural experiments, see for example \cite{LM98}. 
Also due to other recently emerged applications such as check-digit systems \cite{S00,SW10}
and the construction of cryptographic functions \cite{Muratovic2014,SGCGM12},
complete mappings have attracted considerable 
attention, see also \cite{GuCao2015, isikThesis, NidWin2005, TuZeHu2014, Winterhof2014, WuLiHeZh21014, ZhaHuCa2015}.

By a well-known result of Carlitz (1953), all permutation polynomials over~$\Fq$
with $q\ge 3$ 
can be generated by
linear polynomials $ax+b$, $a,b \in \Fq$, $a\neq0$, and {\em inversions} $x^{q-2}=\left\{\begin{array}{cc} 0, &x=0,\\ x^{-1}, & x\ne 0,\end{array}\right.$
see \cite{Carlitz1953} or \cite[Theorem~7.18]{LN97}. Consequently, as pointed out in
\cite{CeMeTo2008}, any permutation $f$ of $\Fq$ can be represented
by a polynomial of the form
\begin{equation}\label{eqn:CarlitzRep}
P_{n}(a_{0},a_{1},...,a_{n+1};x)=(\ldots
((a_0x+a_1)^{q-2}+a_2)^{q-2} \ldots +a_n)^{q-2}+a_{n+1},
\end{equation}
where $a_i \neq 0$, for $i = 0,2,\ldots,n$. Note that this representation is not unique,
and $n$ is not necessarily minimal. Accordingly 
the authors of \cite{aaaw2009} define the {\it Carlitz rank} of a 
permutation polynomial $f$ over $\Fq$ to be the smallest 
integer $n\ge 0$ satisfying $f=P_{n}$ for a permutation $P_{n}$ of the
form (\ref{eqn:CarlitzRep}), and denote it by $Crk(f)$.
In other words, for $q\ge 4$, $Crk(f)=n$ if $f$ is a composition of at least $n$ inversions
$x^{q-2}$ and $n$ or $n+1$ linear polynomials (depending on $a_{n+1}$ being zero or not). This concept, introduced
in the last decade, has already found interesting applications
in diverse areas, see \cite{ACM,AAD,AF}.

The following theorem states the well-known conjecture of Chowla and
Zassenhaus (1968) \cite{Chowla1968}, which was proven by Cohen \cite{Cohen1990} in
1990. 

\textbf{Theorem A.} If $d\geq2$ and $p>(d^2-3d+4)^2$, then there is no
complete mapping polynomial of degree $d$ over $\Fp$.

Note that Cohen's theorem is not true for arbitrary finite fields without further restrictions. For example, for any $0\ne a\in \F_{p^r}$ with 
$a^{(p^r-1)/(p-1)}\ne (-1)^r$ it is easy to see that $ax^p$ is a complete mapping. 


Since the Carlitz rank of a permutation polynomial $f$ over $\Fq$ is an invariant 
of $f$, a natural question to ask is whether a non-existence result, similar to
that stated in Theorem A, can be obtained in terms of the Carlitz rank. 

We define the {\em linearity} ${\cal L}(f)$ of a polynomial $f$ over $\F_q$ by
$${\cal L}(f)=\max_{a,b\in \F_q} |\{ c\in \F_q : f(c)=ac+b\}|.$$
Note that polynomials of large linearity are highly predictable and thus unsuitable in cryptography.

In this paper we show, see Theorem \ref{thm:main} below, that for any 
$n<\lfloor q/2\rfloor$, there is no complete mapping polynomial of Carlitz rank $n$ 
and linearity ${\cal L}(f)<\lfloor (q+5)/2\rfloor$. 

We also answer the following two questions 
that immediately arise. Firstly one wonders how far the non-complete mapping
$f$ in the above setting is from being complete. This question can be quantified
by considering the number $|V_{f+x}|$ of elements in the image of the polynomial
$f+x$. Theorem~\ref{thm2} presents bounds for $|V_{f+x}|$. 
Secondly one would ask if the bound  
$q>2n+1$ can be improved. This is not possible in general, see Example~\ref{ex7} below.  

\section{Preliminaries}

Let $f(x)$ be a permutation polynomial over $\mathbb{F}_{q}$. 
Suppose that $f$ has a representation
$P_{n}$ as in (\ref{eqn:CarlitzRep}) for $n\geq1$. We follow the
notation of \cite{Topuzoglu} and put
$$f(x)=P_{n}(a_{0},a_{1},...,a_{n+1};x).$$

Since we are interested in complete mapping polynomials, the value
of $a_{n+1}$ is irrelevant. Also, by using the substitution
$x\mapsto x-a_{0}^{-1}a_{1}$, we see that the size of the value set
of $f(x)+x$ does not depend on $a_{1}$. Therefore we may
restrict ourselves to the case $a_{1}=a_{n+1}=0$. We relabel the
coefficients accordingly, as $c_{0}=a_{0}$, $c_{i}=a_{i+1}$ for $i=1,..,n-1$, and
 use the notation
\begin {equation}
\label{eqn:f(x)_1}
f(x)=P_{n}(c_{0},...,c_{n-1};x)=:P_n(x).
\end{equation}

The representation of a permutation $f$ as in (\ref {eqn:CarlitzRep}) (or in (\ref{eqn:f(x)_1})) enables approximation of
$f$ by a fractional linear transformation $R_n$ as described below. 

Following the 
terminology of \cite{aaaw2009}, the $n${\it th convergent} $R_{n}(x)$
can be associated to $f$, which is defined as
\begin {equation}
\label{nth Conv} R_{n}(x)=
\frac{\alpha_{n-1}x+\beta_{n-1}}{\alpha_{n}x+\beta_{n}},
\end{equation}
where
\begin{equation*}\label{alpha-beta}
\alpha_{k}=c_{k-1}\alpha_{k-1}+\alpha_{k-2}\hspace{.2in}\mbox{and}\hspace{.2in}
\beta_{k}=c_{k-1}\beta_{k-1}+\beta_{k-2},
\end{equation*}
for $k \geq 2$ and $\alpha_{0}=0,\; \alpha_{1}=c_{0},\; \beta_{0}=1,
\;\beta_{1}=0$.

The set of {\it {poles}} $\mathbf{O}_{n}$ is defined as
\begin{equation*}
\label{poles} 
\mathbf{O_{n}}=\{x_{i}:
x_{i}=\frac{-\beta_{i}}{\alpha_{i}}, \; i=1, \ldots,n \} \subset
\mathbb{F}_{q} \cup \{ \infty \},
\end{equation*}
where the elements of $\mathbf{O_{n}}$ are not necessarily distinct.
We note that
\begin{equation}\label{approx}
f(c)=P_{n}(c)=R_{n}(c)
\quad\mbox{for }c \in\mathbb{F}_{q}\backslash \mathbf{O}_{n}.
\end{equation}

\section{A non-existence result}




In this section we show that any complete mapping must have either high Carlitz rank or high linearity.

\begin{thm}\label{thm:main}
If $f(x)$ is a complete mapping of $\F_q$, then we have either
$${\cal L}(f)\ge \left\lfloor \frac{q+5}{2}\right\rfloor$$
or
$$Crk(f) \ge \left\lfloor \frac{q}{2}\right\rfloor.$$
\end{thm}
\begin{proof}
Let $f(x)$ be of the form $(\ref{eqn:f(x)_1})$ with $n=Crk(f)$ and put $F(x)=f(x)+x$. For $n=0$ we have ${\cal L}(f)=q$. Hence, we may assume $n\ge 1$. 

If $\alpha_n=0$, then $R_n(x)$ defined by $(\ref{nth Conv})$ is a polynomial of degree $1$ with $R_n(c)=f(c)$ for all $c\in \F_q\setminus \On$ by $(\ref{approx})$ and thus
${\cal L}(f)\ge q-n+1$. Since otherwise the result is trivial, we may assume $n\le\lfloor q/2\rfloor-1$ and thus ${\cal L}(f)\ge q+2-\lfloor q/2\rfloor=\lfloor (q+5)/2\rfloor$.

Now we assume $\alpha_n\ne 0$.

We note that the first pole $x_1$ is $0$, since
$\beta_1=0$.
Observe that 
\begin{eqnarray}\label{eqn:first_class_F(c)}
F(c)=R_n(c)+c=\frac{\alpha_{n}{c}^2+(\alpha_{n-1}+\beta_{n})c+\beta_{n-1}}{\alpha_n c +\beta_n}
\end{eqnarray}
for any $c \in \Fq\setminus{\mathbf{O}}_n$. 
It is also easy to show that 

\begin{equation}\label{c_0}
\alpha_{n}\beta_{n-1}-\alpha_{n-1}\beta_{n}=(-1)^{n-1}c_0,
\quad\quad n\ge1.
  \end{equation}

First we assume that $q$ is odd.\\ 
For any $u\in \F_q$ we study the quadratic equation
\begin{equation}\label{quadratic}
R_n(x)+x=u+(\alpha_{n-1}-\beta_n)\alpha_n^{-1}, 
\end{equation}
that is,
\begin{equation}\label{quad}
x^2+(2\alpha_n^{-1}\beta_n-u)x+((-1)^{n-1}c_0+\beta_n^2-u\alpha_n\beta_n)\alpha_n^{-2}=0
\end{equation}
by $(\ref{eqn:first_class_F(c)})$ and $(\ref{c_0})$. 
This equation 
has at most two different solutions $c\in \F_q\setminus\{x_n\}$ and we have exactly two solutions if its discriminant
\begin{equation}\label{discr}
D_u=u^2+4(-1)^nc_0\alpha_n^{-2} 
\end{equation} 
is a square in $\F_q^*$. 
Note that
$$\frac{1+\eta(D_u)}{2}=\left\{\begin{array}{cl} 1, & \mbox{$D_u$ is a square in $\F_q^*$},\\
                                0, & \mbox{$D_u$ is a nonsquare in $\F_q^*$},\\
                                1/2, & D_u=0,
                               \end{array}\right.$$
where $\eta$ is the quadratic character of $\F_q$. 
Moreover, either $D_u=0$ for two values of $u$, that is, $(-1)^{n-1}c_0$ is a square, or there is no value $u$ with $D_u=0$.
Hence, the number $N$ of the elements $u\in \F_q$ for which $D_u$ is a square in $\F_q^*$ can be expressed as
\begin{eqnarray*}
N&=&\frac{1}{2}\sum_{u\in \F_q, D_u\ne 0}(1+\eta(D_u))
=-\frac{1+\eta((-1)^{n-1}c_0)}{2}+\frac{1}{2}\sum_{u\in \F_q} (1+\eta(D_u))\\ 
&=&\frac{q-1-\eta((-1)^{n-1}c_0)}{2}+\frac{1}{2}\sum_{u\in \F_q}\eta(D_u)
=\frac{q-2-\eta((-1)^{n-1}c_0)}{2},
\end{eqnarray*}
by \cite[Theorem~5.48]{LN97}.

Now assume that $F$ is a permutation. Then at least one of these two solutions must be a pole $c\in \On\setminus\{x_n\}$. Hence,  
$$n\ge \frac{q-\eta((-1)^{n-1}c_0)}{2}\ge \frac{q-1}{2}.$$

For even $q$ we can argue similarly. Note that a quadratic equation $x^2+ax+b$ has exactly two solutions whenever $a\ne 0$ and ${\rm Tr}(a^{-2}b)=0$, where
${\rm Tr}$ denotes the absolute trace of $\F_q$, see \cite[Theorem 2.25]{LN97}. We have to determine the number $N$ of $u$ such that $(\ref{quad})$ has two solutions in $\F_q$, that is,
the number of $u\ne 0$ with 
\begin{equation}\label{tr} 0={\rm Tr}\left(\frac{\alpha_n\beta_n u+\beta_n^2+c_0}{\alpha_n^2u^2}\right)
={\rm Tr}\left(\frac{\beta_n}{\alpha_nu}+\frac{\beta_n+c_0^{q/2}}{\alpha_nu}\right)
={\rm Tr}\left(\frac{c_0^{q/2}}{\alpha_nu}\right).
\end{equation}
Since $u\mapsto u^{-1}$ is a bijection of $\F_q^*$ and ${\rm Tr}$ is $2$-to-$1$ on $\F_q$, we get $N=q/2-1$. 
Hence, if $F$ is a permutation, then $\On$ contains at least $n\ge N+1= \frac{q}{2}$ different poles and the result follows.
\end{proof}

Remark. Note that complete mappings of high linearity, that is, polynomials $f(x)$ with $n$th convergent $R_n(x)$ and $\alpha_n=0$ (or $x_n=\infty$) are not suitable for cryptographic applications. 
Hence, in the following we focus on the case $\alpha_n\ne 0$ (or $x_n\ne \infty$).
Note that $\alpha_1\alpha_2\ne 0$ and thus $\infty$ is not a pole if $n=1$ or $n=2$.

Now we provide examples
of complete mappings of Carlitz rank $n=\lfloor q/2\rfloor$ with ${\cal L}(f)<\lfloor (q+5)/2\rfloor$.

\begin{ex}\label{ex7}
 It is easy to check that $f(x)=\gamma (x^4+1)+\gamma^{-1}(x^2+x)\in \F_8[x]$ is a complete mapping of $\F_8=\F_2(\gamma)$, where $\gamma$ is a root of the polynomial $x^3+x+1$ 
 which is irreducible over $\F_2$.
 As a polynomial of degree $4$ its linearity is at most $4$ and by Theorem~\ref{thm:main} its Carlitz rank is at least $4$. Verifying
 $$f(c)= ((((\gamma c)^6+1)^6+\gamma^{-3})^6+1)^6,\quad c\in \F_8,$$
 we see that $Crk(f)=4$ and Theorem~\ref{thm:main} is in general tight in the case of even~$q$.

 Analogously, $f(x)=x^4-x^3+3x^2-x+1\in \F_7[x]$ satisfies
 $$f(c)=(((c^5+3)^5+3)^5,\quad c\in \F_7,$$
 and has Carlitz rank $3$. Hence, the bound of Theorem~\ref{thm:main} is attained for odd $q$, as well.
\end{ex}

Many similar examples lead the authors to believe that there is a complete mapping of $\F_q$ of Carlitz rank $n=\lfloor q/2\rfloor$ and small linearity for infinitely many 
prime powers $q\ge 7$.
This can be checked for $7\le q\le  25$.

\section{The size of $V_{f+x}$ }

In this section we study the set $V_{f+x}=\{f(\delta)+\delta : \delta \in \mathbb{F}_{q}\}$ 
for any $f$ satisfying $(\ref{approx})$ with $\alpha_n\ne 0$. 
Theorem~\ref{thm:main} implies that if $n<\lfloor q/2\rfloor$, we have $|V_{f+x}|<q$. 
Here we aim to determine how large the gap between $q$ and $|V_{f+x}|$ is. 
Theorem~\ref{thm2} below shows that $q-|V_{f+x}| \geq (q-2~ Crk(f)-1)/2$, that is,
it is large if the Carlitz rank of $f$ is small, as one would expect. 
We present the result in a slightly more general form.

\begin{thm}\label{thm2} For $\alpha_{n-1},\beta_{n-1},\alpha_n,\beta_n\in\F_q$ with $\alpha_n\ne 0$ 
and $\alpha_{n-1}\beta_n-\alpha_n\beta_{n-1}\ne 0$, let $F$ be any self-mapping of $\F_q$ satisfying 
\begin{equation}\label{ratio} F(c)=\frac{\alpha_{n-1}c+\beta_{n-1}}{\alpha_nc+\beta_n}+c
\end{equation}
for at least $q-n$ different $c\in \F_q$.
Then we have
\begin{equation*}
\left\lceil\frac{q-n}{2}\right\rceil\leq |V_F|\leq
\min\left\{n+\left\lfloor\frac{q+1}{2}\right\rfloor,q\right\}.
\end{equation*}
\end{thm}

\begin{proof} 
Consider the set $S$ of elements $c\in \F_q$ satisfying $(\ref{ratio})$, which has cardinality $|S|\ge q-n$. At most two different elements of $S$ can have the same value $u$ since $F(c)=u$ is a quadratic equation in $c$ because of the conditions 
on $\alpha_{n-1},\beta_{n-1},\alpha_n,\beta_n$. Therefore,
$|V_F|\ge (q-n)/2$. 
Now the elements of $\F_q\setminus S$ can attain at most $n$ different values of $F$. If $q$ is odd, the discriminant $D_u$ of $F(c)=u$ is a quadratic polynomial in $u$ and is $0$ for 
at most two different values $u\in V_F$. For these two possible $u$ we have exactly one solution $c$ of $F(c)=u$. For all other $u$ we have either two or no solutions.
Hence, the value set of $\frac{\alpha_{n-1}x+\beta_{n-1}}{\alpha_nx+\beta_n}+x$ contains at most $(q+1)/2$ elements and we get
$|V_F|\le n+(q+1)/2$.  
If $q$ is even, the quadratic equation $F(c)=u$ has a unique solution for exactly one $u$ and two or no solutions otherwise. Hence,
we get similarly $|V_F|\le n+q/2$.
\end{proof}

For the special cases $n=1$ and $n=2$ one can provide exact formulas for~$|V_{f+x}|$. 

\begin{prop}
The size of the value set $V_F$ of the polynomial $$F(x)=(c_0x)^{q-2}+x\in \F_q[x],$$ 
$q>2$, with $c_0\ne 0$ is
\begin{equation*}
|V_{F}| = \left\{\begin{array}{cc}(q+1+\eta(c_0)-\eta(-c_0))/2, & \mbox{$q$ odd},\\ q/2, & \mbox{$q$ even},\end{array}\right.
  \end{equation*}
where $\eta$ denotes the quadratic character of $\F_q$. 
\end{prop}

\begin{proof} 
We start with odd $q$.
We have $F(0)=0=F(\pm (-c_0)^{-1/2})$ and thus $F(c)=0$ is attained for $2+\eta(-c_0)$ different $c\in \F_q$.
The discriminant 
$$D_u=u^2-4c_0^{-1}$$ of $x^2-ux+c_0^{-1}$ has no zeros if $c_0$ is a non-square. If $c_0$ is a square, for the two zeros of $D_u$ there is a unique solution $c=u/2$ of $F(c)=u$.
For the remaining $u$ there are two or no solutions of $F(c)=u$. 
Collecting everything we get the result. 

For even $q$ we have $F(0)=F(c_0^{-q/2})=0$ and no further zeros of $F$.
For all $u\ne 0$ there are either two or no solutions of $F(c)=u$ and we get the result.
\end{proof}

\begin{prop}
The size of the value set of 
$F(x)=\big((c_{0}x)^{q-2}+c_{1}\big)^{q-2}+x$, $q>2$, with
$c_{0},c_{1},4c_0+1,c_0+4\ne 0$ is
$$|V_F|=\left\{\begin{array}{ll}\frac{q+2-\eta(4c_0+1)-\eta(c_0^2+4c_0)+\eta(-c_0)}{2}, & c_0\ne -1,\\
                \frac{q-\eta(-3)}{2}, & c_0=-1,
               \end{array}\right.
$$
if $q$ is odd.
For even $q$ and $c_0,c_1\ne 0$, we get
$$|V_F|=\frac{q}{2}+\left\{\begin{array}{ll} Tr(c_0)+Tr(c_0^{-1}),& c_0\ne 1,\\
                Tr(1)-1, &c_0=1,
               \end{array}\right.
$$
where $Tr$ is the absolute trace of $\F_q$ and we identify $\F_2$ with the integers $\{0,1\}$.
\end{prop}
\begin{proof}  Note that 
$\mathbf{O}_{2}=\{0, -(c_{0}c_{1})^{-1}\}$. 
We have $F(0)=c_1^{-1}$ and $$F(-(c_0c_1)^{-1})=-(c_0c_1)^{-1}.$$ 
Note that both values coincide if $c_0=-1$.
$(\ref{quadratic})$ simplifies to $R_2(x)+x=u+c_1^{-1}-(c_0c_1)^{-1}$. Hence, we get $R_2(c)+c=F(0)$ if $u=(c_0c_1)^{-1}=:u_1$ and $R_2(c)+c=F(-(c_0c_1)^{-1})$ if $u=-c_1^{-1}=:u_2$.

Again we deal with odd $q$ first. 

By $(\ref{discr})$ we get the discriminants 
$$D_{u_1}=(4c_0+1)(c_0c_1)^{-2}\quad\mbox{and}\quad D_{u_2}=(c_0+4)c_0(c_0c_1)^{-2}.$$
Hence there are $1+\eta(4c_0+1)$ additional $c$ with $R_2(c)+c=F(0)$ and $1+\eta((c_0+4)c_0)$ additional $c$ with $R_2(c)+c=F(-(c_0c_1)^{-1})$. 
Now verify that there is a $u$, namely $u=(1-c_0)(c_0c_1)^{-1}$, such that $x=0$ is a solution of $(\ref{quad})$. If $c_0=-1$, $x=0$ is the unique solution for this $u$. However, for $x=-(c_0c_1)^{-1}$ there is no such $u$. 
Finally, there are $1+\eta(-c_0)$ values $u$ with $D_u=0$ such that $(\ref{quad})$ has a unique solution. Altogether we have
$$4+\eta(-c_0)+\frac{q-6-\eta(4c_0+1)-\eta((c_0+4)c_0)-\eta(-c_0)}{2}$$
values in $V_F$ if $c_0\ne -1$ and the first result follows.
For $c_0=-1$ we get $|V_F|=2+\frac{q-4-\eta(-3)}{2}$.

Now we consider even $q$.
By $(\ref{tr})$ and 
$$Tr\left(\frac{c_0^{q/2}}{\alpha_2u_1}\right)=Tr(c_0)\quad \mbox{and}\quad
Tr\left(\frac{c_0^{q/2}}{\alpha_2u_2}\right)=Tr(c_0^{-1})$$
the number of $c$ with $F(c)=F(0)$ (including $c=0$) is $3-2Tr(c_0)$ and the number of $c$ with $F(c)=F((c_0c_1)^{-1})$ is
$3-2Tr(c_0^{-1})$. For $u=0$ there is a unique solution $x\ne 0$  of $(\ref{quad})$ if $c_0\ne 1$. Moreover,  $x=0$  is a solution of $(\ref{quad})$ for one $u$ which has already been counted above. 
Hence, we get
$$|V_F|=4+\frac{q-8+2Tr(c_0)+2Tr(c_0^{-1})}{2}$$
if $c_0\ne 1$
and the result follows.

If $c_0=1$ we have $F(0)=F((c_0c_1)^{-1})=c_1^{-1}$ and $c_1^{-1}$ is attained $4-2Tr(c_0)$ times. Moreover, the $u$ with unique solution  $(\ref{quad})$ corresponds to the solution $x=0$.
Hence we get 
$$|V_F|=1+\frac{q-4+2Tr(1)}{2}$$
and the result follows.
\end{proof}

\section{Acknowledgement}

L.I.\ and A.T.\ were supported by TUBITAK project number 114F432.
A.W.\ is partially supported by the Austrian Science Fund FWF Project F5511-N26
which is part of the Special Research Program "Quasi-Monte Carlo Methods: Theory and Applications".


\begin{thebibliography}{999}

\bibitem{aaaw2009} E. Aksoy, A.\c Ce\c smelio\u glu, W. Meidl, A. Topuzo\u glu,
\textit{On the Carlitz rank of a permutation polynomial,}
\emph{Finite Fields and Their Applications} 15 (2009),
428--440.

\bibitem{Carlitz1953} L. Carlitz,
\textit{Permutations in a finite field}, \emph{Proc. American Mathematical
Society} 4 (1953), 538.


\bibitem{Chowla1968} S. Chowla, H. Zassenhaus,
\textit{Some conjectures concerning finite fields}, 
\emph{Norske Videnskabers Selskabs Forhandlinger (Trondheim)} 41 (1968), 34--35.

\bibitem{CeMeTo2008} A. \c Ce\c smelio\u glu, W. Meidl, A. Topuzo\u glu,
\textit{On the cycle structure of permutation polynomials.}
\emph{Finite Fields and Their Applications} 14 (2008),
593--614.

\bibitem{ACM} A. \c Ce\c smelio\u glu, W. Meidl, A. Topuzo\u glu,
\textit{Permutations with prescribed properties},
\emph{Journal of Computational and Applied Mathematics} 259 B (2014), 536--545.

\bibitem{Cohen1990} S.D. Cohen,
\textit{Proof of a conjecture of Chowla and Zassenhaus on
permutation polynomials,} \emph{Canada Mathematical Bulletin} 33
(1990), 230--234.


\bibitem{AAD} D. Gomez-Perez, A. Ostafe, A. Topuzo\u glu,
\textit{On the Carlitz rank of permutations of $\Fq$ and pseudorandom sequences}, 
\emph{Journal of Complexity} 30 (2014), 279--289.

\bibitem{GuCao2015} X. Guangkui, X. Cao, \textit{Complete
permutation polynomials over finite fields of odd characteristic,}
\emph{Finite Fields and Their Applications} 31 (2015),
228--240.

\bibitem{isikThesis} L. I\c s\i k, \textit{On complete mappings and value sets of polynomials over finite fields},
 PhD Thesis. Sabanc{\i} University, 2015.

\bibitem{LM98} C.F. Laywine, G. Mullen, 
Discrete mathematics using Latin squares. 
Wiley-Interscience Series in Discrete Mathematics and Optimization. A Wiley-Interscience Publication. John Wiley \& Sons, Inc., New York, 1998.
 
\bibitem{LN97} R. Lidl, H. Niederreiter, Finite fields. Second edition. Encyclopedia of Mathematics and its Applications, 20. Cambridge University Press, Cambridge, 1997.
 
\bibitem{Mann1942} H.B. Mann,
\textit{The construction of orthogonal Latin squares}, \emph{Annals
of Mathematical Statistics} 13 (1942), 418--423.


\bibitem{Muratovic2014} A. Muratovic-Ribic, E. Pasalic,
\textit{ A note on complete mapping polynomials over finite fields
and their applications in cryptography,} \emph{Finite Fields and
Their Applications} 25 (2014), 306--315.

\bibitem{NR} H. Niederreiter, K.H. Robinson,
\textit{Complete mappings of finite fields,} \emph{Journal of Australian Mathematical
Society} A 33 (1982), 197--212.

\bibitem{NidWin2005} H. Niederreiter, A. Winterhof,
\textit{Cyclotomic $\mathcal{R}$-orthomorphisms of finite fields},
\emph{Discrete Mathematics} 295 (2005), 161--171.

\bibitem{AF} F. Pausinger, A. Topuzo\u glu,
\textit{Permutations of finite fields and uniform distribution modulo 1}, in
H. Niederrreiter, A. Ostafe, D. Panario, A. Winterhof (eds.), \emph{Algebraic Curves and Finite Fields}, Radon Series on Applied and Computational Mathematics 16 (2014), 145--160.

\bibitem{S00} R.-H. Schulz, \textit{On check digit systems using anti-symmetric mappings}, in Numbers, information and complexity (Bielefeld, 1998), 295--310, Kluwer Acad. Publ., Boston, MA, 2000.

\bibitem{SW10} R. Shaheen, A. Winterhof, \textit{Permutations of finite fields for check digit systems}, \emph{Des. Codes Cryptogr.} 57 (2010), 361--371.

\bibitem{SGCGM12} P. St\u anic\u a, S. Gangopadhyay, A. Chaturvedi, A.K. Gangopadhyay, S. Maitra,
\textit{Investigations on bent and negabent functions via the nega-Hadamard transform}, \emph{IEEE Trans. Inf. Theory} 58 (2012), 4064--4072.

\bibitem{Topuzoglu} A. Topuzo\u glu, \textit{Carlitz rank of permutations of
finite fields: A survey}, \emph{Journal of Symbolic Computation}  64
(2014), 53--66.

\bibitem{TuZeHu2014} Z. Tu, X. Zeng, L. Hu,
\textit{Several classes of complete permutation polynomials,}
\emph{Finite Fields and Their Applications} 25 (2014),
182--193.

\bibitem{Winterhof2014} A. Winterhof,
\textit{Generalizations of complete mappings of finite
fields and some applications}, \emph{ Journal of Symbolic
Computation} 64 (2014), 42--52.

\bibitem{WuLiHeZh21014} G. Wu, N. Li, T. Helleseth, Y. Zhang,
\textit{Some classes of monomial complete permutation polynomials
over finite fields of characteristic two,} \emph{Finite Fields and
Their Applications} 28 (2014), 148--165.

\bibitem{ZhaHuCa2015} Z. Zha, L. Hu, X. Cao,
\textit{Constructing permutations and complete permutations over
finite fields via subfield-valued polynomials}, \emph{Finite Fields
and Their Applications} 31 (2015), 162--177.


\end{thebibliography}
\end{document}